\documentclass[11pt]{article}
\usepackage{graphicx}
\usepackage{amsfonts}
\usepackage{theorem}

\newtheorem{Lem}{Lemma}[section]
\newtheorem{Def}[Lem]{Definition}
\newtheorem{The}[Lem]{Theorem}
\newtheorem{Prop}[Lem]{Proposition}
\newtheorem{Cor}[Lem]{Corollary}

\newtheorem{Rem}[Lem]{Remark}

\input{tcilatex}

\begin{document}

\title{Some inequalities on generalized entropies}
\author{S. Furuichi$^{1}$\thanks{%
E-mail:furuichi@chs.nihon-u.ac.jp}, N. Minculete$^{2}$\thanks{%
E-mail:minculeten@yahoo.com} and F.-C. Mitroi$^{3}$\thanks{%
E-mail:fcmitroi@yahoo.com} \\
$^1${\small Department of Computer Science and System Analysis,}\\
{\small College of Humanities and Sciences, Nihon University,}\\
{\small 3-25-40, Sakurajyousui, Setagaya-ku, Tokyo, 156-8550, Japan}\\
$^2${\small \lq\lq Dimitrie Cantemir\rq\rq University, Bra\c{s}ov, 500068,
Rom{a}nia}\\
$^3${\small University of Craiova, Department of Mathematics,}\\
{\small Street A. I. Cuza 13, Craiova, RO-200585, Romania}}
\date{}
\maketitle

\textbf{Abstract.} We give several inequalities on generalized entropies
involving Tsallis entropies, using some inequalities obtained by
improvements of Young's inequality. We also give a generalized Han's
inequality. \vspace{3mm}

 \vspace{3mm}

 \vspace{%
3mm}


\section{Introduction}

We start from the weighted quasilinear mean for some continuous and strictly
monotonic function $\psi :I\rightarrow \mathbb{R}$, defined by 
\begin{equation}
M_{\psi }(x_{1},x_{2},\cdots ,x_{n})\equiv \psi ^{-1}\left(
\sum_{j=1}^{n}p_{j}\psi (x_{j})\right) ,
\end{equation}%
where $\sum_{j=1}^{n}p_{j}=1$, $p_{j}>0$, $x_{j}\in I$ for $j=1,2,\cdots ,n$
and $n\in \mathbb{N}$. If we take $\psi (x)=x$, then $M_{\psi
}(x_{1},x_{2},\cdots ,x_{n})$ coincides with the weighted arithmetic mean $%
A(x_{1},x_{2},\cdots ,x_{n})\equiv \sum_{j=1}^{n}p_{j}x_{j}$. If we also
take $\psi (x)=\log (x)$, then $M_{\psi }(x_{1},x_{2},\cdots ,x_{n})$
coincides with the weighted geometric mean $G(x_{1},x_{2},\cdots
,x_{n})\equiv \prod_{j=1}^{n}x_{j}^{p_{j}}$.

If $\psi (x)=x$ and $x_{j}=\ln _{q}\frac{1}{p_{j}}$, then $M_{\psi
}(x_{1},x_{2},\cdots ,x_{n})$ is equal to Tsallis entropy \cite{Tsa}: 
\begin{equation}
H_{q}(p_{1},p_{2},\cdots ,p_{n})\equiv -\sum_{j=1}^{n}p_{j}^{q}\ln
_{q}p_{j}=\sum_{j=1}^{n}p_{j}\ln _{q}\frac{1}{p_{j}},\,\,(q\geq 0,q\neq 1)
\end{equation}%
where $\{p_{1},p_{2},\cdots ,p_{n}\}$ is a probability distribution with $%
p_{j}>0$ for all $j=1,2,\cdots ,n$ and the $q-$logarithmic function for $x>0$
is defined by $\ln _{q}(x)\equiv \frac{x^{1-q}-1}{1-q}$ which uniformly
converges to the usual logarithmic function $\log (x)$ in the limit $%
q\rightarrow 1$. Therefore Tsallis entropy conveges to Shannon entropy in
the limit $q\rightarrow 1$: 
\begin{equation}
\lim_{q\rightarrow 1}H_{q}(p_{1},p_{2},\cdots
,p_{n})=H_{1}(p_{1},p_{2},\cdots ,p_{n})\equiv -\sum_{j=1}^{n}p_{j}\log
p_{j}.
\end{equation}%
Thus we find that the Tsallis entropy is one of the generalizations of
Shannon entropy. It is known that the R\'{e}nyi entropy \cite{Ren} is also a
generalization of Shannon entropy. Here, we review the quasilinear entropy 
\cite{AD} as another generalization of Shannon entropy. For a continuous and
strictly monotonic function $\phi $ on $(0,1]$, the quasilinear entropy is
given by 
\begin{equation}
I^{\phi }(p_{1},p_{2},\cdots ,p_{n})\equiv -\log \phi ^{-1}\left(
\sum_{j=1}^{n}p_{j}\phi (p_{j})\right) .  \label{sec1_04}
\end{equation}%
If we take $\phi (x)=\log \left( x\right) $ in (\ref{sec1_04}), then we have 
$I^{\log }(p_{1},p_{2},\cdots ,p_{n})=H_{1}(p_{1},p_{2},\cdots ,p_{n})$. We
may redefine quasilinear entropy by 
\begin{equation}
I_{1}^{\psi }(p_{1},p_{2},\cdots ,p_{n})\equiv \log \psi ^{-1}\left(
\sum_{j=1}^{n}p_{j}\psi \left( \frac{1}{p_{j}}\right) \right) ,
\label{sec1_05}
\end{equation}%
for a continuous and strictly monotonic function $\psi $ on $(0,\infty )$.
If we take $\psi (x)=\log \left( x\right) $ in (\ref{sec1_05}), we have $%
I_{1}^{\log }(p_{1},p_{2},\cdots ,p_{n})=H_{1}(p_{1},p_{2},\cdots ,p_{n})$.
The case $\psi (x)=x^{1-q}$ is also useful in practice,$\ $since we
recapture R\'{e}nyi entropy, namely $I_{1}^{x^{1-q}}(p_{1},p_{2},\cdots
,p_{n})=R_{q}(p_{1},p_{2},\cdots ,p_{n})$ where R\'{e}nyi entropy \cite{Ren}
is defined by 
\begin{equation}
R_{q}(p_{1},p_{2},\cdots ,p_{n})\equiv \frac{1}{1-q}\log \left(
\sum_{j=1}^{n}p_{j}^{q}\right) .  \label{sec2_eq15}
\end{equation}

\begin{Def}
For a continuous and strictly monotonic function $\psi $ on $(0,\infty)$and
two probability distributions $\{p_{1},p_{2},\cdots ,p_{n}\}$ and $%
\{r_{1},r_{2},\cdots ,r_{n}\}$ with $p_{j}>0,r_{j}>0$ for all $j=1,2,\cdots
,n$, the quasilinear relative entropy is defined by%
\begin{equation}
D_{1}^{\psi }(p_{1},p_2,\cdots ,p_{n}||r_{1},r_2,\cdots ,r_{n})\equiv -\log
\psi ^{-1}\left( \sum_{j=1}^{n}p_{j}\psi \left( \frac{r_{j}}{p_{j}}\right)
\right) .
\end{equation}
\end{Def}

The quasilinear relative entropy coincides to Shannon relative entropy if $%
\psi (x)=\log \left( x\right) ,$ i.e. 
\[
D_{1}^{\log }(p_{1},p_{2},\cdots ,p_{n}||r_{1},r_{2},\cdots
,r_{n})=-\sum_{j=1}^{n}p_{j}\log \frac{r_{j}}{p_{j}}=D_{1}(p_{1},p_{2},%
\cdots ,p_{n}||r_{1},r_{2},\cdots ,r_{n}).
\]%
We denote by $R_{q}(p_{1},p_{2},\cdots ,p_{n}||r_{1},r_{2},\cdots ,r_{n})$
the R\'{e}nyi relative entropy \cite{Ren} defined by 
\begin{equation}
R_{q}(p_{1},p_{2},\cdots ,p_{n}||r_{1},r_{2},\cdots ,r_{n})\equiv \frac{1}{%
q-1}\log \left( \sum_{j=1}^{n}p_{j}^{q}r_{j}^{1-q}\right) .
\end{equation}%
This is another particular case of quasilinear relative entropy, namely for $%
\psi (x)=x^{1-q}\ $we have 
\begin{eqnarray*}
D_{1}^{x^{1-q}}(p_{1},p_{2},\cdots ,p_{n}||r_{1},r_{2},\cdots ,r_{n})
&=&-\log \left( \sum_{j=1}^{n}p_{j}\left( \frac{r_{j}}{p_{j}}\right)
^{1-q}\right) ^{\frac{1}{1-q}}=\frac{1}{q-1}\log \left(
\sum_{j=1}^{n}p_{j}^{q}r_{j}^{1-q}\right) \newline
\\
&=&R_{q}(p_{1},p_{2},\cdots ,p_{n}||r_{1},r_{2},\cdots ,r_{n}).
\end{eqnarray*}%
We denote by 
\begin{equation}
D_{q}(p_{1},p_{2},\cdots ,p_{n}||r_{1},r_{2},\cdots ,r_{n})\equiv
\sum_{j=1}^{n}p_{j}^{q}(\ln _{q}p_{j}-\ln _{q}r_{j})=-\sum_{j=1}^{n}p_{j}\ln
_{q}\frac{r_{j}}{p_{j}}
\end{equation}%
the Tsallis relative entropy. Tsallis relative entropy conveges to the usual
relative entropy (divergence, K-L information) in the limit $q\rightarrow 1$%
: 
\begin{eqnarray}
\lim_{q\rightarrow 1}D_{q}(p_{1},p_{2},\cdots ,p_{n}||r_{1},r_{2},\cdots
,r_{n}) &=&D_{1}(p_{1},p_{2},\cdots ,p_{n}||r_{1},r_{2},\cdots ,r_{n}) 
\nonumber \\
&\equiv &\sum_{j=1}^{n}p_{j}(\log p_{j}-\log r_{j}).
\end{eqnarray}%
See \cite%
{Tsa1,Tsa2,Tsa3,Furu_2004,Furu_2005,Furu_2006,Furu_2007,Furu_2008,Furu_2009,Furu_2010,Furu_2011}
and references therein for recent advances and applications on the Tsallis
entropy. We easily find that the Tsallis relative entropy is a special case
of Csisz\'{a}r $f$-divergence \cite{Csi1,Csi2,Csi3} defined for a convex
function $f$ on $(0,\infty )$ with $f(1)=0$ by 
\begin{equation}
D_{f}(p_{1},p_{2},\cdots ,p_{n}||r_{1},r_{2},\cdots ,r_{n})\equiv
\sum_{j=1}^{n}r_{j}f\left( \frac{p_{j}}{r_{j}}\right) ,
\end{equation}%
since $f(x)=-x\ln _{q}\left( 1/x\right) $ is convex on $(0,\infty )$,
vanishes at $x=1$ and 
\[
D_{-x\ln _{q}(1/x)}(p_{1},p_{2},\cdots ,p_{n}||r_{1},r_{2},\cdots
,r_{n})=D_{q}(p_{1},p_{2},\cdots ,p_{n}||r_{1},r_{2},\cdots ,r_{n}).
\]%
Furthermore, we define the dual function with respect to a convex function $f
$ by 
\begin{equation}
f^{\ast }(t)=tf\left( \frac{1}{t}\right) 
\end{equation}%
for $t>0$. Then the function $f^{\ast }(t)$ is also convex on $(0,\infty )$.
In addition, we define the $f$-divergence for \textit{incomplete}
probability distributions $\{a_{1},a_{2},\cdots ,a_{n}\}$ and $%
\{b_{1},b_{2}\cdots ,b_{n}\}$ where $a_{i}>0$ and $b_{i}>0$, in the
following way: 
\begin{equation}
\widetilde{D_{f^{\ast }}}(a_{1},a_{2},\cdots ,a_{n}||b_{1},b_{2},\cdots
,b_{n})\equiv \sum_{j=1}^{n}a_{j}f^{\ast }\left( \frac{b_{j}}{a_{j}}\right) .
\end{equation}

On the other hand, the studies on refinements for Young's inequality have
given a great progress in the papers \cite%
{BK,Dra,KM,Ald2,Ald1,Mit,Furu1,Furu2,Min1,Min2,FM,MF}. In the present paper,
we give some inequalities on the Tsallis entropies applying two type
inequalities obtained in \cite{Mit,Min1}. In addition, we give the
generalized Han's inequality for the Tsallis entropy in the final section.


\section{Tsallis quasilinear entropy and Tsallis quasilinear relative entropy%
}

As an analogy with (\ref{sec1_05}), we may define the following entropy.

\begin{Def}
For a continuous and strictly monotonic function $\psi $ on $(0,\infty)$ and 
$q \geq 0$ with $q \neq 1$, Tsallis quasilinear entropy ($q$-quasilinear
entropy) is defined by 
\begin{equation}
I_q^{\psi }(p_{1},p_{2},\cdots ,p_{n}) \equiv \ln _{q}\psi ^{-1}\left(
\sum_{j=1}^{n}p_{j}\psi \left(\frac{1}{p_{j}}\right)\right),
\end{equation}
where $\{p_1,p_2,\cdots,p_n\}$ is a probability distribution with $p_j>0$
for all $j=1,2,\cdots,n$.
\end{Def}

We notice that if $\psi $ does not depend on $q$ then $\lim_{q\rightarrow
1}I_{q}^{\psi }(p_{1},p_{2},\cdots ,p_{n})=I_{1}^{\psi }(p_{1},p_{2},\cdots
,p_{n}).$

For $x>0$ and $q\geq 0$ with $q\neq 1$, we define the $q$-exponential
function as the inverse function of the $q$-logarithmic function by $\exp
_{q}(x)\equiv \left\{ 1+(1-q)x\right\} ^{1/(1-q)}$, if $1+(1-q)x>0$,
otherwise it is undefined. 
If we take $\psi (x)=\ln _{q}(x)$ then we have $I_{q}^{\ln
_{q}}(p_{1},p_{2},\cdots ,p_{n})=H_{q}(p_{1},p_{2},\cdots ,p_{n}).$
Furthermore, we have 
\begin{eqnarray*}
I_{q}^{x^{1-q}}(p_{1},p_{2},\cdots ,p_{n}) &=&\ln _{q}\left(
\sum_{j=1}^{n}p_{j}p_{j}^{q-1}\right) ^{\frac{1}{1-q}}=\ln _{q}\left(
\sum_{j=1}^{n}p_{j}^{q}\right) ^{\frac{1}{1-q}} \\
&=&\frac{\left[ \left( \sum_{j=1}^{n}p_{j}^{q}\right) ^{\frac{1}{1-q}}\right]
^{1-q}-1}{1-q}=\frac{\sum_{j=1}^{n}\left( p_{j}^{q}-p_{j}\right) }{1-q}%
=H_{q}(p_{1},p_{2},\cdots ,p_{n}).
\end{eqnarray*}

\begin{Prop}
Tsallis quasilinear entropy is nonnegative: 
\[
I_q^{\psi }(p_{1},p_{2},\cdots ,p_{n}) \geq 0. 
\]
\end{Prop}

\textit{Proof}: We assume that $\psi $ is an increasing function. Then we
have $\psi \left( \frac{1}{p_{j}}\right) \geq \psi (1)$ from $\frac{1}{p_{j}}%
\geq 1$ for $p_{j}>0$ for all $j=1,2,\cdots ,n$. Thus we have $%
\sum_{j=1}^{n}p_{j}\psi \left( \frac{1}{p_{j}}\right) \geq \psi (1)$ which
implies $\psi ^{-1}\left( \sum_{j=1}^{n}p_{j}\psi \left( \frac{1}{p_{j}}%
\right) \right) \geq 1$, since $\psi ^{-1}$ is also increasing. For the case
that $\psi $ is a decreasing function, we can prove it similarly.

\hfill \hbox{\rule{6pt}{6pt}}

We note here that the $q$-exponential function gives us the following
connection between R\'{e}nyi entropy and Tsallis entropy \cite{Mas}: 
\begin{equation}
\exp R_{q}(p_{1},p_{2},\cdots ,p_{n})=\exp _{q}H_{q}(p_{1},p_{2},\cdots
,p_{n}).  \label{sec2_eq14}
\end{equation}%
We should note here $\exp _{q}H_{q}(p_{1},p_{2},\cdots ,p_{n})$ is always
defined, since we have 
\[
1+(1-q)H_{q}(p_{1},p_{2},\cdots ,p_{n})=\sum_{j=1}^{n}p_{j}^{q}>0.
\]%
From (\ref{sec2_eq14}), we have the following proposition.

\begin{Prop}
\label{sec2_prop1} Let $\mathcal{A}\equiv \left\{ \mathcal{A}%
_{i}:i=1,2,\cdots ,k\right\} $ be a partition of $\{1,2,\cdots ,n\}$ and put 
$p_{i}^{\mathcal{A}}\equiv \sum_{j\in \mathcal{A}_{i}}p_{j}$. Then we have 
\begin{eqnarray}
&&\sum_{j=1}^{n}p_{j}^{q}\geq \sum_{j=1}^{k}\left( p_{j}^{\mathcal{A}%
}\right) ^{q},\,\,\,(0\leq q\leq 1), \\
&&\sum_{j=1}^{n}p_{j}^{q}\leq \sum_{j=1}^{k}\left( p_{j}^{\mathcal{A}%
}\right) ^{q},\,\,\,(1\leq q).
\end{eqnarray}
\end{Prop}

\textit{Proof}: We use the generalized Shannon additivity (which is often
called $q$-additivity) for Tsallis entropy (see \cite{Furu_2005} for
example): 
\begin{equation}
\hspace*{-8mm}H_{q}(x_{11},\cdots ,x_{nm_{n}})=H_{q}(x_{1},\cdots
,x_{n})+\sum_{i=1}^{n}x_{i}^{q}H_{q}\left( \frac{x_{i1}}{x_{i}},\cdots ,%
\frac{x_{im_{i}}}{x_{i}}\right) .  \label{sec2_prop1_eq01}
\end{equation}%
where $x_{ij}\geq 0$, $x_{i}=\sum_{j=1}^{m_{i}}x_{ij},\,\,(i=1,\cdots
,n;j=1,\cdots ,m_{i})$. Thus we have 
\begin{equation}
H_{q}(p_{1},p_{2},\cdots ,p_{n})\geq H_{q}\left( p_{1}^{\mathcal{A}},p_{2}^{%
\mathcal{A}},\cdots ,p_{k}^{\mathcal{A}}\right) ,  \label{sec2_prop1_eq02}
\end{equation}%
since the second term of the right hand side in (\ref{sec2_prop1_eq01}) is
nonnegative, because of the nonnegativity of Tsallis entropy. Thus we have 
\begin{eqnarray*}
\exp R_{q}(p_{1},p_{2},\cdots ,p_{n}) &=&\exp _{q}H_{q}(p_{1},p_{2},\cdots
,p_{n}) \\
&\geq &\exp _{q}H_{q}\left( p_{1}^{\mathcal{A}},p_{2}^{\mathcal{A}},\cdots
,p_{k}^{\mathcal{A}}\right) \\
&=&\exp R_{q}\left( p_{1}^{\mathcal{A}},p_{2}^{\mathcal{A}},\cdots ,p_{k}^{%
\mathcal{A}}\right) ,
\end{eqnarray*}%
since $\exp _{q}$ is a monotone increasing function. Hence the inequality 
\begin{equation}
R_{q}(p_{1},p_{2},\cdots ,p_{n})\geq R_{q}\left( p_{1}^{\mathcal{A}},p_{2}^{%
\mathcal{A}},\cdots ,p_{k}^{\mathcal{A}}\right) ,  \label{sec2_prop1_eq03}
\end{equation}%
holds, which proves the present proposition.

\hfill \hbox{\rule{6pt}{6pt}}

\begin{Def}
For a continuous and strictly monotonic function $\psi $ on $(0,\infty)$ and
two probability distributions $\{p_{1},p_{2},\cdots ,p_{n}\}$ and $%
\{r_{1},r_{2},\cdots ,r_{n}\}$ with $p_{j}>0,r_{j}>0$ for all $j=1,2,\cdots
,n$, the Tsallis quasilinear relative entropy is defined by 
\begin{equation}
D_{q}^{\psi }(p_{1},p_2,\cdots ,p_{n}||r_{1},r_2,\cdots ,r_{n}) \equiv -\ln
_{q}\psi ^{-1}\left( \sum_{j=1}^{n}p_{j}\psi \left( \frac{r_{j}}{p_{j}}%
\right) \right) .
\end{equation}
\end{Def}

For $\psi (x)=\ln _{q}\left( x\right) $ the Tsallis quasilinear relative
entropy becomes Tsallis relative entropy, that is 
\[
D_{q}^{\ln _{q}}(p_{1},p_{2},\cdots ,p_{n}||r_{1},r_{2},\cdots
,r_{n})=-\sum_{j=1}^{n}p_{j}\ln _{q}\frac{r_{j}}{p_{j}}=D_{q}(p_{1},p_{2},%
\cdots ,p_{n}||r_{1},r_{2},\cdots ,r_{n}),
\]%
and for $\psi (x)=x^{1-q}$, we have 
\begin{eqnarray}
&&D_{q}^{x^{1-q}}(p_{1},p_{2},\cdots ,p_{n}||r_{1},r_{2},\cdots ,r_{n})=-\ln
_{q}\left( \sum_{j=1}^{n}p_{j}\left( \frac{r_{j}}{p_{j}}\right)
^{1-q}\right) ^{\frac{1}{1-q}}=-\ln _{q}\left(
\sum_{j=1}^{n}p_{j}^{q}r_{j}^{1-q}\right) ^{\frac{1}{1-q}}  \nonumber \\
&=&\frac{-\left\{ \left[ \left( \sum_{j=1}^{n}p_{j}^{q}r_{j}^{1-q}\right) ^{%
\frac{1}{1-q}}\right] ^{1-q}-1\right\} }{1-q}=\frac{\sum_{j=1}^{n}\left(
p_{j}-p_{j}^{q}r_{j}^{1-q}\right) }{1-q}  \nonumber \\
&=&D_{q}(p_{1},p_{2},\cdots ,p_{n}||r_{1},r_{2},\cdots ,r_{n}).
\end{eqnarray}%
We give a sufficient condition on nonnegativity of Tsallis quasilinear
relative entropy.

\begin{Prop}
If $\psi$ is a concave increasing function or a convex decreasing function,
then we have nonnegativity of Tsallis quasilinear relative entropy: 
\[
D_{q}^{\psi }(p_{1},p_2,\cdots ,p_{n}||r_{1},r_2,\cdots ,r_{n}) \geq 0. 
\]
\end{Prop}

\textit{Proof}: We firstly assume that $\psi$ is a concave increasing
function. The concavity of $\psi$ shows that we have $\psi\left(%
\sum_{j=1}^np_j\frac{r_j}{p_j} \right) \geq \sum_{j=1}^n p_j \psi\left(\frac{%
r_j}{p_j}\right)$ which is equivalent to $\psi (1) \geq \sum_{j=1}^n p_j
\psi\left(\frac{r_j}{p_j}\right).$ From the assumption, $\psi^{-1}$ is also
increasing so that we have $1 \geq \psi^{-1}\left( \sum_{j=1}^n p_j
\psi\left(\frac{r_j}{p_j}\right) \right). $ Therefore we have $-\ln_q
\psi^{-1}\left( \sum_{j=1}^n p_j \psi\left(\frac{r_j}{p_j}\right) \right)
\geq 0, $ since $\ln_q x$ is increasing and $\ln_q(1)=0$. For the case that $%
\psi$ is a convex decreasing function, we can prove similarly the
nonnegativity of Tsallis quasilinear relative entropy.

\hfill \hbox{\rule{6pt}{6pt}}

\begin{Rem}
The following two functions satisfy the sufficient condition in the above
proposition.

\begin{itemize}
\item[(i)] $\psi(x) = \ln_q x$ for $q \geq 0, q\neq 1$.

\item[(ii)] $\psi(x)=x^{1-q}$ for $q \geq 0, q\neq 1$.
\end{itemize}
\end{Rem}

It is notable that the following identity holds 
\begin{equation}
\exp R_{q}(p_{1},p_{2},\cdots ,p_{n}||r_{1},r_{2},\cdots ,r_{n})=
\exp_{2-q}D_{q}(p_{1},p_{2},\cdots ,p_{n}||r_{1},r_{2},\cdots ,r_{n}).
\label{sec2_eq16}
\end{equation}
We should note here $\exp_{2-q} D_{q}(p_{1},p_{2},\cdots
,p_{n}||r_{1},r_{2},\cdots ,r_{n}) $ is always defined, since we have 
\[
1+(q-1) D_{q}(p_{1},p_{2},\cdots ,p_{n}||r_{1},r_{2},\cdots
,r_{n})=\sum_{j=1}^n p_j^q r_j^{1-q} >0. 
\]
We also find that (\ref{sec2_eq16}) implies the monotonicity of R\'enyi
relative entropy.

\begin{Prop}
\label{sec2_prop2} Under the same assumptions with Proposition \ref%
{sec2_prop1} and $r_i^{\mathcal{A}} \equiv \sum_{j\in\mathcal{A}_i} r_j$, we
have 
\begin{equation}
R_{q}(p_{1},p_{2},\cdots ,p_{n}||r_{1},r_{2},\cdots ,r_{n}) \geq
R_q\left(p_1^{\mathcal{A}},p_2^{\mathcal{A}},\cdots,p_k^{\mathcal{A}}||r_1^{%
\mathcal{A}},r_2^{\mathcal{A}},\cdots,r_k^{\mathcal{A}}\right).
\end{equation}
\end{Prop}

\textit{Proof}: We recall that Tsallis relative entropy is a special case of 
$f$-divergence so that it has same properties with $f$-divergence. Since $%
\exp _{2-q}$ is a monotone increasing function for $0\leq q \leq 2$ and $f$%
-divergence has a monotonicity \cite{Csi1,Csi3}, we have 
\begin{eqnarray*}
\exp R_{q}(p_{1},p_{2},\cdots ,p_{n}||r_{1},r_{2},\cdots ,r_{n}) &=&\exp
_{2-q}D_{q}(p_{1},p_{2},\cdots ,p_{n}||r_{1},r_{2},\cdots ,r_{n}) \\
&\geq &\exp _{2-q}D_{q}\left( p_{1}^{\mathcal{A}},p_{2}^{\mathcal{A}},\cdots
,p_{k}^{\mathcal{A}}||r_{1}^{\mathcal{A}},r_{2}^{\mathcal{A}},\cdots ,r_{k}^{%
\mathcal{A}}\right) \\
&=&\exp R_{q}\left( p_{1}^{\mathcal{A}},p_{2}^{\mathcal{A}},\cdots ,p_{k}^{%
\mathcal{A}}||r_{1}^{\mathcal{A}},r_{2}^{\mathcal{A}},\cdots ,r_{k}^{%
\mathcal{A}}\right) ,
\end{eqnarray*}%
which proves the statement.

\hfill \hbox{\rule{6pt}{6pt}}


\section{Inequalities for Tsallis quasilinear entropy and $f$-divergence}

In this section, we give inequalities for Tsallis quasilinear entropy and $f$%
-divergence. For this purpose, we review the results obtained in \cite{Mit}
as one of generalizations of refined Young's inequality.

\begin{Prop}
\textbf{(\cite{Mit})} \label{prop_01} For two probability vectors $\mathbf{p}%
=\{p_{1},p_{2},\cdots ,p_{n}\}$ and $\mathbf{r}=\{r_{1},r_{2},\cdots
,r_{n}\} $ such that $p_{j}>0$, $r_{j}>0$, $\sum_{j=1}^{n}p_{j}=%
\sum_{j=1}^{n}r_{j}=1$ and $\mathbf{x}=\{x_{1},x_{2},\cdots ,x_{n}\}$ such
that $x_{i}\geq 0$, we have 
\begin{equation}
\min_{1\leq i\leq n}\left\{ \frac{r_{i}}{p_{i}}\right\} T(f,\mathbf{x},%
\mathbf{p})\leq T(f,\mathbf{x},\mathbf{r})\leq \max_{1\leq i\leq n}\left\{ 
\frac{r_{i}}{p_{i}}\right\} T(f,\mathbf{x},\mathbf{p}),
\end{equation}%
where 
\begin{equation}
T(f,\mathbf{x},\mathbf{p})\equiv \sum_{j=1}^{n}p_{j}f(x_{j})-f\left( \psi
^{-1}\left( \sum_{j=1}^{n}p_{j}\psi (x_{j})\right) \right) ,
\end{equation}%
for continuous increasing function $\psi :I\rightarrow I$ and a function $%
f:I\rightarrow J$ such that 
\begin{equation}
f(\psi ^{-1}((1-\lambda )\psi (a)+\lambda \psi (b)))\leq (1-\lambda
)f(a)+\lambda f(b)
\end{equation}%
for any $a,b\in I$ and any $\lambda \in \lbrack 0,1]$.
\end{Prop}

We have the following inequalities on Tsallis quasilinear entropy and
Tsallis entropy.

\begin{The}
\label{sec2_the1} For $q\geq 0$, a continuous and strictly monotonic
function $\psi $ on $(0,\infty )$ and a probability distribution $%
\{r_{1},r_{2},\cdots ,r_{n}\}$ with $r_{j}>0$ for all $j=1,2,\cdots ,n$, we
have 
\begin{eqnarray*}
0 &\leq &n\min_{1\leq i\leq n}\{r_{i}\}\left\{ \ln _{q}\left( \psi
^{-1}\left( \frac{1}{n}\sum_{j=1}^{n}\psi \left( \frac{1}{r_{j}}\right)
\right) \right) -\frac{1}{n}\sum_{j=1}^{n}\ln _{q}\frac{1}{r_{j}}\right\} \\
&\leq &I_{q}^{\psi }(r_{1},r_{2},\cdots ,r_{n})-H_{q}\left(
r_{1},r_{2},\cdots ,r_{n}\right) \\
&\leq &n\max_{1\leq i\leq n}\{r_{i}\}\left\{ \ln _{q}\left( \psi ^{-1}\left( 
\frac{1}{n}\sum_{j=1}^{n}\psi \left( \frac{1}{r_{j}}\right) \right) \right) -%
\frac{1}{n}\sum_{j=1}^{n}\ln _{q}\frac{1}{r_{j}}\right\} \\
&&
\end{eqnarray*}
\end{The}

\textit{Proof}: If we take the uniform distribution $\mathbf{p}=\left\{\frac{%
1}{n},\cdots,\frac{1}{n}\right\}\equiv \mathbf{u}$ in Proposition \ref%
{prop_01}, then we have 
\begin{equation}  \label{sec2_the1_01}
n \min_{1\leq i \leq n} \left\{ r_i \right\} T_n(f,\mathbf{x},\mathbf{u})
\leq T_n(f,\mathbf{x},\mathbf{r}) \leq n \max_{1\leq i \leq n}\left\{r_i
\right\} T_n(f,\mathbf{x},\mathbf{u}),
\end{equation}
(which coincides with Theorem 3.3 in \cite{Mit}). In the inequalities (\ref%
{sec2_the1_01}), we put $f(x)=-\ln_q(x)$ and $x_j=\frac{1}{r_j}$ for any $%
j=1,2,\cdots,n$, then we obtain the statement.

\hfill \hbox{\rule{6pt}{6pt}}

\begin{Cor}
\label{sec2_cor1} For $q\geq 0$ and a probability distribution $%
\{r_{1},r_{2},\cdots ,r_{n}\}$ with $r_{j}>0$ for all $j=1,2,\cdots ,n$, we
have 
\begin{eqnarray}
&&\hspace*{-1cm}0\leq n\min_{1\leq i\leq n}\{r_{i}\}\left\{ \ln _{q}\left( 
\frac{1}{n}\sum_{j=1}^{n}\frac{1}{r_{j}}\right) -\frac{1}{n}%
\sum_{j=1}^{n}\ln _{q}\frac{1}{r_{j}}\right\} \leq \ln
_{q}n-H_{q}(r_{1},r_{2},\cdots ,r_{n})  \nonumber \\
&&\hspace*{4cm}\leq n\max_{1\leq i\leq n}\{r_{i}\}\left\{ \ln _{q}\left( 
\frac{1}{n}\sum_{j=1}^{n}\frac{1}{r_{j}}\right) -\frac{1}{n}%
\sum_{j=1}^{n}\ln _{q}\frac{1}{r_{j}}\right\} ,  \label{sec2_cor1_01}
\end{eqnarray}
\end{Cor}

\textit{Proof}: Put $\psi(x)=x$ in Theorem \ref{sec2_the1}.

\hfill \hbox{\rule{6pt}{6pt}}

\begin{Rem}
Corollary \ref{sec2_cor1} improves the well-known inequalities $0 \leq
H_q(r_1,r_2,\cdots,r_n) \leq \ln_q n$. If we take the limit $q \to 1$, the
inequalities (\ref{sec2_cor1_01}) recover Proposition 1 in \cite{Dra}.
\end{Rem}

We also have the following inequalities.

\begin{The}
\label{sec2_the2} For two probability distributions $\mathbf{p}=\left\{
p_{1},p_{2},\cdots ,p_{n}\right\} $ and $\mathbf{r}=\{r_{1},r_{2},\cdots
,r_{n}\}$, and an \textit{incomplete} probability distribution $\mathbf{t}%
=\{t_{1},t_{2},\cdots ,t_{n}\}$ with $t_{j}\equiv \frac{p_{j}^{2}}{r_{j}}$,
we have 
\begin{eqnarray}
&&0\leq \min_{1\leq i\leq n}\left\{ \frac{r_{i}}{p_{i}}\right\} \left( 
\widetilde{D_{f^{\ast }}}(\mathbf{t}||\mathbf{p})-f\left(
\sum_{j=1}^{n}t_{j}\right) \right)   \nonumber \\
&\leq &D_{f}(\mathbf{p}||\mathbf{r})\leq \max_{1\leq i\leq n}\left\{ \frac{%
r_{i}}{p_{i}}\right\} \left( \widetilde{D_{f^{\ast }}}(\mathbf{t}||\mathbf{p}%
)-f\left( \sum_{j=1}^{n}t_{j}\right) \right) .  \label{sec2_the2_ineq30}
\end{eqnarray}
\end{The}

\textit{Proof}: Put $x_{j}=\frac{p_{j}}{r_{j}}$ in Proposition \ref{prop_01}
with $\psi (x)=x$. Since we have the relation 
\[
\sum_{j=1}^{n}p_{j}f\left( \frac{p_{j}}{r_{j}}\right) =\sum_{j=1}^{n}p_{j}%
\frac{p_{j}}{r_{j}}f^{\ast }\left( \frac{r_{j}}{p_{j}}\right)
=\sum_{j=1}^{n}t_{j}f^{\ast }\left( \frac{p_{j}}{t_{j}}\right) , 
\]%
we have the statement.

\hfill \hbox{\rule{6pt}{6pt}}

\begin{Cor}
\textbf{(\cite{Dra})} Under the same assumption as in Theorem \ref{sec2_the2}%
, we have 
\begin{eqnarray*}
&&0\leq \min_{1\leq i\leq n}\left\{ \frac{r_{i}}{p_{i}}\right\} \left( \log
\left( \sum_{j=1}^{n}t_{j}\right) -D_{1}(\mathbf{p}||\mathbf{r})\right)  \\
&\leq &D_{1}(\mathbf{r}||\mathbf{p})\leq \max_{1\leq i\leq n}\left\{ \frac{%
r_{i}}{p_{i}}\right\} \left( \log \left( \sum_{j=1}^{n}t_{j}\right) -D_{1}(%
\mathbf{p}||\mathbf{r})\right) .
\end{eqnarray*}
\end{Cor}

\textit{Proof}: If we take $f(x)=-\log \left( x\right) $ in Theorem \ref%
{sec2_the2}, then we have 
\[
D_{f}(\mathbf{p}||\mathbf{r})=-\sum_{j=1}^{n}r_{j}\log \frac{p_{j}}{r_{j}}%
=\sum_{j=1}^{n}r_{j}\log \frac{r_{j}}{p_{j}}=D_{1}(\mathbf{r}||\mathbf{p}).
\]%
Since $f^{\ast }(x)=x\log \left( x\right) $ and $t_{j}=\frac{p_{j}^{2}}{r_{j}%
}$, we also have 
\begin{eqnarray*}
&&\widetilde{D_{f^{\ast }}}(\mathbf{t}||\mathbf{p})-f\left(
\sum_{j=1}^{n}t_{j}\right) =\sum_{j=1}^{n}t_{j}\frac{p_{j}}{t_{j}}\log \frac{%
p_{j}}{t_{j}}+\log \left( \sum_{j=1}^{n}t_{j}\right)
=\sum_{j=1}^{n}p_{j}\log \frac{r_{j}}{p_{j}}+\log \left(
\sum_{j=1}^{n}t_{j}\right)  \\
&=&-\sum_{j=1}^{n}p_{j}\log \frac{p_{j}}{r_{j}}+\log \left(
\sum_{j=1}^{n}t_{j}\right) =\log \left( \sum_{j=1}^{n}t_{j}\right) -D_{1}(%
\mathbf{p}||\mathbf{r}).
\end{eqnarray*}

\hfill \hbox{\rule{6pt}{6pt}}



\section{Inequalities for Tsallis entropy}

We firstly give Lagrange's identity \cite{Wei}, to establish an alternative
generalization of refined Young's inequality.

\begin{Lem}
\textbf{(Lagrange's identity)} \label{sec3_lem1} For two vectors $%
\{a_1,a_2,\cdots,a_n\}$ and $\{b_1,b_2,\cdots,b_n\}$, we have 
\begin{eqnarray}
\left(\sum_{k=1}^n a_k^2\right)\left(\sum_{k=1}^n
b_k^2\right)-\left(\sum_{k=1}^n a_k b_k\right)^2 &=& \frac{1}{2}%
\sum_{i=1}^n\sum_{j=1}^n \left(a_ib_j-a_jb_i\right)^2  \nonumber \\
&=& \sum_{1\leq i<j \leq n} \left(a_ib_j-a_jb_i\right)^2.
\end{eqnarray}
\end{Lem}

\begin{The}
\label{sec3_the1} Let $f:I \to \mathbb{R}$ be a twice differentiable
function such that there exist real constants $m$ and $M$ so that $0\leq
m\leq f^{\prime\prime}(x)\leq M$ for any $x\in I$. Then we have 
\begin{eqnarray}
\frac{m}{2} \sum_{1\leq i<j \leq n} p_ip_j\left( x_j- x_i \right)^2 &\leq&
\sum_{j=1}^n p_jf(x_j) -f\left(\sum_{j=1}^np_j x_j\right)  \nonumber \\
&\leq & \frac{M}{2}\sum_{1\leq i<j \leq n} p_ip_j\left( x_j- x_i \right)^2
\label{sec3_the1_01}
\end{eqnarray}
where $p_j >0$ with $\sum_{j=1}^n p_j =1$ and $x_j\in I$ for all $%
j=1,2,\cdots,n$.
\end{The}

\textit{Proof}: We consider the function $g:I\rightarrow \mathbb{R}$ defined
by $g(x)\equiv f(x)-\frac{m}{2}x^{2}$. Since we have $g^{\prime \prime
}(x)=f^{\prime \prime }(x)-m\geq 0$, $g$ is a convex function. Applying
Jensen's inequality, we thus have 
\begin{equation}
\sum_{j=1}^{n}p_{j}g(x_{j})\geq g\left( \sum_{j=1}^{n}p_{j}x_{j}\right)
\label{sec3_the1_02}
\end{equation}%
where $p_{j}>0$ with $\sum_{j=1}^{n}p_{j}=1$ and $x_{j}\in I$ for all $%
j=1,2,\cdots ,n$. From the inequality (\ref{sec3_the1_02}), we have 
\begin{eqnarray*}
\sum_{j=1}^{n}p_{j}f(x_{j})-f\left( \sum_{j=1}^{n}p_{j}x_{j}\right) &\geq &%
\frac{m}{2}\left\{ \sum_{j=1}^{n}p_{j}x_{j}^{2}-\left(
\sum_{j=1}^{n}p_{j}x_{j}\right) ^{2}\right\} \\
&=&\frac{m}{2}\left\{ \left( \sum_{j=1}^{n}p_{j}\right) \left(
\sum_{j=1}^{n}p_{j}x_{j}^{2}\right) -\left( \sum_{j=1}^{n}p_{j}x_{j}\right)
^{2}\right\} \\
&=&\frac{m}{2}\sum_{1\leq i<j\leq n}\left( \sqrt{p_{i}}\sqrt{p_{j}}x_{j}-%
\sqrt{p_{j}}\sqrt{p_{i}}x_{i}\right) ^{2} \\
&=&\frac{m}{2}\sum_{1\leq i<j\leq n}p_{i}p_{j}\left( x_{j}-x_{i}\right) ^{2}.
\end{eqnarray*}%
In the above calculations, we used Lemma \ref{sec3_lem1}. Thus we proved the
first part of the inequalities. Similarly, one can prove the second part of
the inequalities, putting the function $h:I\rightarrow \mathbb{R}$ defined
by $h(x)\equiv \frac{M}{2}x^{2}-f(x)$. We omit the details.

\hfill \hbox{\rule{6pt}{6pt}}

\begin{Lem}
\label{sec3_lem2} For $\{p_1,p_2,\cdots,p_n\} $ with $p_j >0$ and $%
\sum_{j=1}^n p_j =1$, and $\{x_1,x_2,\cdots,x_n\}$ with $x_j>0$, we have 
\begin{equation}
\sum_{1\leq i<j \leq n} p_ip_j\left( x_j- x_i \right)^2 =
\sum_{j=1}^np_j\left(x_j-\sum_{i=1}^np_ix_i\right)^2.
\end{equation}
\end{Lem}

\textit{Proof}: We denote 
\[
\bar{x}=\sum_{i=1}^{n}p_{i}x_{i}. 
\]%
The left side term becomes 
\begin{eqnarray*}
\sum_{1\leq i<j\leq n}p_{i}p_{j}\left( x_{j}-x_{i}\right) ^{2} &=&\frac{1}{2}%
\sum_{i=1}^{n}\sum_{j=1}^{n}p_{i}p_{j}\left( x_{j}-x_{i}\right) ^{2}=\frac{1%
}{2}\sum_{i=1}^{n}\sum_{j=1}^{n}p_{i}p_{j}\left(
x_{j}^{2}+x_{i}^{2}-2x_{j}x_{i}\right) \\
&=&\frac{1}{2}\sum_{i=1}^{n}\sum_{j=1}^{n}p_{i}p_{j}x_{j}^{2}+\frac{1}{2}%
\sum_{i=1}^{n}\sum_{j=1}^{n}p_{i}p_{j}x_{i}^{2}-\sum_{i=1}^{n}%
\sum_{j=1}^{n}p_{i}p_{j}x_{j}x_{i} \\
&=&\frac{1}{2}\sum_{i=1}^{n}p_{i}\sum_{j=1}^{n}p_{j}x_{j}^{2}+\frac{1}{2}%
\sum_{i=1}^{n}p_{i}x_{i}^{2}\sum_{j=1}^{n}p_{j}-\sum_{i=1}^{n}p_{i}x_{i}%
\sum_{j=1}^{n}p_{j}x_{j} \\
&=&\sum_{j=1}^{n}p_{j}x_{j}^{2}-\bar{x}^{2}.
\end{eqnarray*}%
Similarly, a straightforward computation yields%
\begin{eqnarray*}
\sum_{j=1}^{n}p_{j}\left( x_{j}-\sum_{i=1}^{n}p_{i}x_{i}\right) ^{2}
&=&\sum_{j=1}^{n}p_{j}\left( x_{j}^{2}-2x_{j}\bar{x}+\bar{x}^{2}\right)
=\sum_{j=1}^{n}p_{j}x_{j}^{2}-2\bar{x}^{2}+\bar{x}^{2} \\
&=&\sum_{j=1}^{n}p_{j}x_{j}^{2}-\bar{x}^{2}.
\end{eqnarray*}%
This concludes the proof.

\hfill \hbox{\rule{6pt}{6pt}}

\begin{Cor}
\label{sec3_cor0} Under the assumptions of Theorem \ref{sec3_the1}, we have 
\begin{eqnarray}
\frac{m}{2}\sum_{j=1}^{n}p_{j}\left( x_{j}-\sum_{i=1}^{n}p_{i}x_{i}\right)
^{2} &\leq &\sum_{j=1}^{n}p_{j}f(x_{j})-f\left(
\sum_{j=1}^{n}p_{j}x_{j}\right)  \nonumber \\
&\leq &\frac{M}{2}\sum_{j=1}^{n}p_{j}\left(
x_{j}-\sum_{i=1}^{n}p_{i}x_{i}\right) ^{2}.
\end{eqnarray}
\end{Cor}

\begin{Rem}
Corollary \ref{sec3_cor0} gives a similar form with Cartwright-Field's
inequality \cite{CF}: 
\begin{eqnarray}
\frac{1}{2M^{\prime}} \sum_{j=1}^n p_j\left(x_j-\sum_{i=1}^n p_i
x_i\right)^2 &\leq& \sum_{j=1}^n p_j x_j-\prod_{j=1}^n x_j^{p_j}  \nonumber
\\
&\leq& \frac{1}{2m^{\prime}} \sum_{j=1}^n p_j\left(x_j-\sum_{i=1}^n p_i
x_i\right)^2  \label{sec3_rem0_01}
\end{eqnarray}
where $p_j>0$ for all $j=1,2,\cdots,n$ and $\sum_{j=1}^n p_j=1$, $%
m^{\prime}\equiv \min\{x_1,x_2,\cdots,x_n\}>0$ and $M^{\prime}\equiv\max%
\{x_1,x_2,\cdots,x_n\}$.
\end{Rem}

We also have the following inequalities for Tsallis entropy.

\begin{The}
\label{sec3_the2} For two probability distributions $\{p_{1},p_{2},\cdots
,p_{n}\}$ and $\{r_{1},r_{2},\cdots ,r_{n}\}$ with $p_{j}>0$, $r_{j}>0$ and $%
\sum_{j=1}^{n}p_{j}=\sum_{j=1}^{n}r_{j}=1$, we have 
\begin{eqnarray}
&&\ln _{q}\left( \sum_{j=1}^{n}\frac{p_{j}}{r_{j}}\right) -\ln _{q}n+\frac{%
m_q}{2}\sum_{1\leq i<j\leq n}p_{i}p_{j}\left( \frac{1}{p_{j}}-\frac{1}{p_{i}}%
\right) ^{2}-\frac{M_q}{2}\sum_{1\leq i<j\leq n}p_{i}p_{j}\left( \frac{1}{%
r_{j}}-\frac{1}{r_{i}}\right) ^{2}  \nonumber \\
&\leq &\sum_{j=1}^{n}p_{j}\ln _{q}\frac{1}{r_{j}}-\sum_{j=1}^{n}p_{j}\ln _{q}%
\frac{1}{p_{j}}  \nonumber \\
&\leq &\ln _{q}\left( \sum_{j=1}^{n}\frac{p_{j}}{r_{j}}\right) -\ln _{q}n+%
\frac{M_q}{2}\sum_{1\leq i<j\leq n}p_{i}p_{j}\left( \frac{1}{p_{j}}-\frac{1}{%
p_{i}}\right) ^{2}-\frac{m_q}{2}\sum_{1\leq i<j\leq n}p_{i}p_{j}\left( \frac{%
1}{r_{j}}-\frac{1}{r_{i}}\right) ^{2},  \nonumber \\
&&
\end{eqnarray}%
where $m_q$ and $M_q$ are positive numbers depending on the parameter $q\geq
0$ and satisfying $m_q\leq q r_{j}^{-q-1} \leq M_q$ and $m_q \leq q
p_{j}^{-q-1} \leq M_q$ for all $j=1,2,\cdots ,n$.
\end{The}

\textit{Proof}: Applying Theorem \ref{sec3_the1} for the convex function $%
-\ln _{q}(x)$ and $x_{j}=\frac{1}{r_{j}}$, we have 
\begin{eqnarray}
\frac{m_q}{2}\sum_{1\leq i<j\leq n}p_{i}p_{j}\left( \frac{1}{r_{j}}-\frac{1}{%
r_{i}}\right) ^{2} &\leq &-\sum_{j=1}^{n}p_{j}\ln _{q}\frac{1}{r_{j}}+\ln
_{q}\left( \sum_{j=1}^{n}\frac{p_{j}}{r_{j}}\right)  \nonumber \\
&\leq &\frac{M_q}{2}\sum_{1\leq i<j\leq n}p_{i}p_{j}\left( \frac{1}{r_{j}}-%
\frac{1}{r_{i}}\right) ^{2},  \label{sec3_the2_01}
\end{eqnarray}%
since the second derivative of $-\ln _{q}(x)$ is $qx^{-q-1}$. Putting $%
r_{j}=p_{j}$ for all $j=1,2,\cdots ,n$ in the inequalities (\ref%
{sec3_the2_01}), it follows 
\begin{eqnarray}
\frac{m_q}{2}\sum_{1\leq i<j\leq n}p_{i}p_{j}\left( \frac{1}{p_{j}}-\frac{1}{%
p_{i}}\right) ^{2} &\leq &-\sum_{j=1}^{n}p_{j}\ln _{q}\frac{1}{p_{j}}+\ln
_{q}n  \nonumber \\
&\leq &\frac{M_q}{2}\sum_{1\leq i<j\leq n}p_{i}p_{j}\left( \frac{1}{p_{j}}-%
\frac{1}{p_{i}}\right) ^{2}.  \label{sec3_the2_02}
\end{eqnarray}%
From the inequalities (\ref{sec3_the2_01}) and (\ref{sec3_the2_02}), we have
the statement.

\hfill \hbox{\rule{6pt}{6pt}}

\begin{Rem}
The first part of the inequalities (\ref{sec3_the2_02}) gives another
improvement of the well-known inequalities $0\leq
H_{q}(r_1,r_2,\cdots,r_n)\leq \ln _{q}n$.
\end{Rem}

\begin{Cor}
For two probability distributions $\{p_{1},p_{2},\cdots ,p_{n}\}$ and $%
\{r_{1},r_{2},\cdots ,r_{n}\}$ with $p_{j}>0$, $r_{j}>0$ and $%
\sum_{j=1}^{n}p_{j}=\sum_{j=1}^{n}r_{j}=1$, we have 
\begin{eqnarray}
&&\log \left( \sum_{j=1}^{n}\frac{p_{j}}{r_{j}}\right) -\log n+\frac{m_1}{2}%
\sum_{1\leq i<j\leq n}p_{i}p_{j}\left( \frac{1}{p_{j}}-\frac{1}{p_{i}}%
\right) ^{2}-\frac{M_1}{2}\sum_{1\leq i<j\leq n}p_{i}p_{j}\left( \frac{1}{%
r_{j}}-\frac{1}{r_{i}}\right) ^{2}  \nonumber \\
&\leq &\sum_{j=1}^{n}p_{j}\log \frac{1}{r_{j}}-\sum_{j=1}^{n}p_{j}\log \frac{%
1}{p_{j}}  \nonumber \\
&\leq &\log \left( \sum_{j=1}^{n}\frac{p_{j}}{r_{j}}\right) -\log n+\frac{M_1%
}{2}\sum_{1\leq i<j\leq n}p_{i}p_{j}\left( \frac{1}{p_{j}}-\frac{1}{p_{i}}%
\right) ^{2}-\frac{m_1}{2}\sum_{1\leq i<j\leq n}p_{i}p_{j}\left( \frac{1}{%
r_{j}}-\frac{1}{r_{i}}\right) ^{2},  \nonumber \\
&&  \label{4}
\end{eqnarray}
where $m_1$ and $M_1$ are positive numbers satisfying $m_1\leq r_{j}^{-2}
\leq M_1$ and $m_1 \leq p_{j}^{-2} \leq M_1$ for all $j=1,2,\cdots ,n$.
\end{Cor}

\textit{Proof}: Take the limit $q \to 1$ in Theorem \ref{sec3_the2}.

\hfill \hbox{\rule{6pt}{6pt}}

\begin{Rem}
The second part of the inequalities (\ref{4}) gives the reverse inequality
for the so-called information inequality \cite[Theorem 2.6.3]{Cov}: 
\begin{equation}
0\leq \sum_{j=1}^{n}p_{j}\log \frac{1}{r_{j}}-\sum_{j=1}^{n}p_{j}\log \frac{1%
}{p_{j}}  \label{sec3_rem2_01}
\end{equation}%
which is equivalent to the non-negativity of the relative entropy: 
\[
D_{1}(p_{1},p_{2},\cdots ,p_{n}||r_{1},r_{2},\cdots ,r_{n})\geq 0. 
\]
\end{Rem}

Using the inequality (\ref{sec3_rem2_01}), we derive the following result.

\begin{Prop}
For two probability distributions $\{p_{1},p_{2},\cdots ,p_{n}\}$ and $%
\{r_{1},r_{2},\cdots ,r_{n}\}$ with $0<p_{j}<1$, $0<r_{j}<1$ and $%
\sum_{j=1}^{n}p_{j}=\sum_{j=1}^{n}r_{j}=1$, we have 
\begin{equation}
\sum_{j=1}^{n}\left( 1-p_{j}\right) \log \frac{1}{1-p_{j}}\leq
\sum_{j=1}^{n}\left( 1-p_{j}\right) \log \frac{1}{1-r_{j}}.
\end{equation}
\end{Prop}

\textit{Proof}: In the inequality (\ref{sec3_rem2_01}), we put $p_{j}=\frac{%
1-p_{j}}{n-1}$ and $r_{j}=\frac{1-r_{j}}{n-1}$ which satisfy $\sum_{j=1}^{n}%
\frac{1-p_{j}}{n-1}=\sum_{j=1}^{n}\frac{1-r_{j}}{n-1}=1$. Then we have the
present proposition.

\hfill \hbox{\rule{6pt}{6pt}}


\section{A generalized Han's inequality}

In order to state our result, we give the definitions of the Tsallis
conditional entropy and the Tsallis joint entropy.

\begin{Def}
\textbf{(\cite{Dar,Furu_2006})} For the conditional probability $p(x_i\vert
y_j)$ and the joint probability $p(x_i,y_j)$, we define the Tsallis
conditional entropy and the Tsallis joint entropy by 
\begin{equation}
H_q(\mathbf{x}\vert \mathbf{y}) \equiv -\sum_{i,j} p(x_i,y_j)^q \ln_q
p(x_i\vert y_j), \quad (q \geq 0, q \neq 1),  \label{jc}
\end{equation}
and 
\begin{equation}
H_q(\mathbf{x},\mathbf{y}) \equiv -\sum_{i,j} p(x_i,y_j)^q \ln_q p(x_i,y_j),
\quad (q \geq 0, q \neq 1).  \label{jc2}
\end{equation}
\end{Def}

We summarize briefly the following chain rules representing relations
between Tsallis conditional entropy and Tsallis joint entropy.

\begin{Prop}
\label{prop1} \textbf{(\cite{Dar,Furu_2006})}Assume that ${\mathbf{x},%
\mathbf{y}}$ are probability distributions. Then 
\begin{equation}
H_{q}\left( {\mathbf{x},\mathbf{y}}\right) =H_{q}\left( \mathbf{x}\right)
+H_{q}\left( {\mathbf{y}\left\vert \mathbf{x}\right. }\right) .
\label{chain1}
\end{equation}
\end{Prop}

Proposition \ref{prop1} implied the following propositions.

\begin{Prop}
\label{prop2}\textbf{(\cite{Furu_2006})} Suppose $\mathbf{x}_{1},\mathbf{x}%
_{2},\cdots ,\mathbf{x}_{n}$ are probability distributions. Then 
\begin{equation}
H_{q}\left( {\mathbf{x}_{1},\mathbf{x}_{2},\cdots ,\mathbf{x}_{n}}\right)
=\sum\limits_{i=1}^{n}{H_{q}\left( {\mathbf{x}_{i}\left\vert {\mathbf{x}%
_{i-1},\cdots ,\mathbf{x}_{1}}\right. }\right) }.  \label{chain_gen}
\end{equation}
\end{Prop}

\begin{Prop}
\label{prop3} \textbf{(\cite{Dar,Furu_2006})} For $q\geq 1$, two probability
distributions $\mathbf{x}$ and $\mathbf{y}$, we have the following
inequality: 
\begin{equation}
H_{q}\left( {\mathbf{x}|\mathbf{y}}\right) \leq H_{q}\left( \mathbf{x}%
\right) .  \label{sub1}
\end{equation}
\end{Prop}

Consequently we have the following self-bounding property of Tsallis joint
entropy.

\begin{The}
\textbf{(Generalized Han's inequality)}\label{the1} Let $\mathbf{x}_1,%
\mathbf{x}_2,\cdots ,\mathbf{x}_n$ be probability distributions. Then for $q
\geq 1$, we have the following inequality: 
\[
H_q(\mathbf{x}_1,\cdots,\mathbf{x}_n)\leq \frac{1}{n-1}\sum_{i=1}^n H_q(%
\mathbf{x}_1,\cdots,\mathbf{x}_{i-1},\mathbf{x}_{i+1},\cdots,\mathbf{x}_n). 
\]
\end{The}

\textit{Proof}: Since the Tsallis joint entropy has a symmetry: $H_{q}(%
\mathbf{x},\mathbf{y})=H_{q}(\mathbf{y},\mathbf{x})$, we have 
\begin{eqnarray*}
H_{q}(\mathbf{x}_{1},\cdots ,\mathbf{x}_{n}) &=&H_{q}(\mathbf{x}_{1},\cdots ,%
\mathbf{x}_{i-1},\mathbf{x}_{i+1},\cdots ,\mathbf{x}_{n})+H_{q}(\mathbf{x}%
_{i}|\mathbf{x}_{1},\cdots ,\mathbf{x}_{i-1},\mathbf{x}_{i+1},\cdots ,%
\mathbf{x}_{n}) \\
&\leq &H_{q}(\mathbf{x}_{1},\cdots ,\mathbf{x}_{i-1},\mathbf{x}_{i+1},\cdots
,\mathbf{x}_{n})+H_{q}(\mathbf{x}_{i}|\mathbf{x}_{1},\cdots ,\mathbf{x}%
_{i-1}),
\end{eqnarray*}%
by the use of Proposition \ref{prop1} and Proposition \ref{prop3}. Summing
both sides on $i$ from $1$ to $n$, we have 
\begin{eqnarray*}
nH_{q}(\mathbf{x}_{1},\cdots ,\mathbf{x}_{n}) &=&\sum_{i=1}^{n}H_{q}(\mathbf{%
x}_{1},\cdots ,\mathbf{x}_{i-1},\mathbf{x}_{i+1},\cdots ,\mathbf{x}%
_{n})+\sum_{i=1}^{n}H_{q}(\mathbf{x}_{i}|\mathbf{x}_{1},\cdots ,\mathbf{x}%
_{i-1},\mathbf{x}_{i+1},\cdots ,\mathbf{x}_{n}) \\
&\leq &\sum_{i=1}^{n}H_{q}(\mathbf{x}_{1},\cdots ,\mathbf{x}_{i-1},\mathbf{x}%
_{i+1},\cdots ,\mathbf{x}_{n})+H_{q}(\mathbf{x}_{1},\cdots ,\mathbf{x}_{n}),
\end{eqnarray*}%
due to Proposition \ref{prop2}. Therefore we have the present proposition.
\hfill \hbox{\rule{6pt}{6pt}}

\begin{Rem}
Theorem \ref{the1} recovers the original Han's inequality \cite{Han,BLB}, if
we take the limit as $q \to 1$.
\end{Rem}


\section*{Acknowledgements}

The author (S.F.) was supported in part by the Japanese Ministry of
Education, Science, Sports and Culture, Grant-in-Aid for Encouragement of
Young Scientists (B), 20740067. The author (N.M.) was supported in part by
the Romanian Ministry of Education, Research and Innovation through the PNII
Idei project 842/2008. The author (F.-C. M.) was supported by CNCSIS Grant $%
420/2008.$

\end{document}